\def\intav{-\hskip -1.1em\int}
\input amstex
\magnification=\magstep1
\documentstyle{amsppt}
\NoRunningHeads
\vsize=8.75 true in

\parindent=.3in

\topmatter
\title
A generalization of Riesz's uniqueness theorem 
\endtitle

\author
Enrique Villamor\\
Department of Mathematics\\
Florida International University\\
Miami, FL 33199
\\
\endauthor

\address
Dept. of Mathematics, Florida International University, Miami, FL 
33199.
\endaddress

\email 
villamor{\@}solix.fiu.edu
\endemail

\keywords
Capacity
\endkeywords

\subjclass
  30C65 
\endsubjclass

\abstract
There have been, over the last 8 years, a number of far reaching of the famous original F. 
and M. Riesz's uniqueness theorem that states that if a bounded analytic function in the unit disc 
of the complex plane $\Bbb C$ has the same radial limit in a set of positive Lebesgue measure 
on its boundary, then the function has to be constant. First Beurling [B], considering the case of 
non-constant meromorphic functions mapping the unit disc on a Riemann surface of finite spherical 
area, was able to prove that if such a function showed an appropriate behavior in the neighborhood 
of the limit value where the function maps a set on the boundary of the unit disc, then those sets 
have capacity zero. Here the capacity considered is the logarithmic linear capacity. 
The author of the present note in [V], was able to weakened beurling condition on the limit 
value. Later Jenkins in [J], showed that in the presence of such a local condition on the limiting 
value, the global behavior of Riemann surface is irrelevant and at the same time he gave an improved 
and sharper condition.\par
Those results where quite restrictive in a two folded way, namely, they were in dimension $n=2$ and 
the regularity requirements on the treated functions were quite strong, analyticity and meromorphicity. 
Koskela in [K], was able to remove those two restrictions by proving a uniqueness result for functions 
in $ACL^p(\Bbb B^n)$ for values of $p$ in the interval $(1,n]$ and satisfying a condition on 
the limit value very similar in nature to the one of Jenkins in dimension 2. In particular, Koskela's 
result recovers Jenkins in the case $p=n=2$. He proves that a continuous function in the Sobolev space 
$W^{1,p}(\Bbb B^n)$ (here $\Bbb B^n$ is the unit ball of $\Bbb R^n$ and $1<p\leq n$) vanishes identically 
provided
$$\int_{\|u(x) -a\| <\epsilon} \|\nabla u(x) \|^p\,\,dx=O(\epsilon^{p}\,\,(\log({1\over\epsilon}))^{p-1})$$
as $\epsilon\to 0$ and there is a set $E$ on $\partial \Bbb B^n$ of positive $p$-capacity such that 
each $x\in E$ is a terminal point of some rectifiable curve along which the function $u$ tends to $a$.\par
	
Koskela also shows in his paper that this result is sharp in the sense that $(\log({1\over\epsilon}))^{p-1}$ 
can not be replaced by $(\log({1\over\epsilon}))^{p-1+\delta}$ for any positive $\delta$ (even if $u$ is 
assumed to be continuous in the closure of $\Bbb B^n$).\par

\endabstract

\endtopmatter

\document

\beginsection{$\S$1. Introduction.}

Mizuta in [M] showed that under the same hypothesis on the function $u$, if 
$$\int_{\|u(x) -a\| <\epsilon} \|\nabla u(x) \|^p\,\,dx=O(\epsilon^{p}\,\,\phi(\epsilon))$$
as $\epsilon\to 0$, where $\phi$ is a positive nonincreasing function on the interval $(0,\infty)$ 
satisfying the following conditions
$$(1)\,\,\,\,\,A^{-1}\,\,\phi(r)\leq \phi(r^2)\leq A\,\, \phi(r)$$
for every $r>0$ and $A$ a positive constant and 
$$(2)\,\,\,\,\,\int_0^1 [\phi(r)]^{{1\over{1-p}}}\,\,\,r^{-1}\,\,dr=\infty,$$
then if there is a set $E$ on $\partial \Bbb B^n$ of 
positive $p$-capacity such that 
each $x\in E$ is a terminal point of some rectifiable curve along which the function $u$ tends to $a$; the 
function $u$ vanishes identically on $\Bbb B^n$. It is easy to observe that the function $\phi(\epsilon)= 
(\log({1\over\epsilon}))^{p-1}$ satisfies the two conditions in [M].\par
	
Last, Miklyukov and Vuorinen in [MV] showed that if we define
$$I(\epsilon)= \int_{\|u(x) -a\| <\epsilon} \|\nabla u(x) \|^p\,\,dx.$$
If the integral $I(\epsilon)$ satisfies one of the conditions
$$(3)\,\,\,\,\,\int_0 ({1\over{I'(\epsilon)}})^{{1\over{p-1}}}\,\,d\epsilon =\infty;$$
or
$$(4)\,\,\,\,\,\int_0 ({\epsilon \over{I(\epsilon)}})^{{1\over{p-1}}}\,\,d\epsilon =\infty;$$
or there exists a nonnegative function $f(\epsilon)$ satisfying conditions
$$(5)\,\,\,\,\,I(\epsilon)\leq \epsilon^p\,\,(f(\epsilon))^{p-1},$$
for every $0<\epsilon<{1\over 2}$; and
$$(6)\,\,\,\,\,\sum_{k=0}^\infty {1\over{f(2^{-k})}}=\infty;$$
or
$$(7)\,\,\,\,\,\liminf_{\epsilon\to 0} {{I(\epsilon)}\over {\epsilon^p}}<\infty;$$
then again the function $u$ is identically equal to 0.\par
	
Again, it is not difficult to show that this result generalizes the one in [M].\par
	
In this paper we are going to give a more general condition on the integral 
$\int_{\|u(x) -a\| <\epsilon} \|\nabla u(x) \|^p\,\,dx$ under which we will still be able to deduce that 
the function $u$ vanishes identically and this condition is more general than the ones appearing in 
[M] and [MV]. Let us remark here that the results of Koskela, Mizuta and Miklyukov and Vuorinen 
are for real valued functions, while our result is going to be more in the spirit of the initial results of 
M. and F. Riesz, Beurling and Jenkins, where they considered functions from the complex plane into the 
complex plane. \par
In that spirit, our results will hold for functions defined in the unit ball of $\Bbb R^n$ into $\Bbb R^n$. Later, 
we will see how this general result includes the results in [K], [M] and [MV].
In this section we will introduce several definitions that will be needed in the rest of the paper. 
Let us start by recalling the definition of monotone function (in 
this paper we 
consider only continuous monotone functions).
\proclaim{Definition 1.1}
	Let $\Omega\subset \Bbb R^n$ be an open set. A continuous function 
$u\colon \Omega\to \Bbb R$ is monotone 
(in the sense of Lebesgue) if 
$$\max_{\bar D} u(x)=\max_{\partial D} u(x)$$
and 
$$\min_{\bar D}u(x)=\min_{\partial D} u(x)$$ 
hold whenever $D$ is a domain with compact closure $\bar D\subset 
\Omega$.
\endproclaim

The Sobolev space $W^{1,p}(\Bbb B^n)$ is defined in [HKM, Chapter 
1]. It consists of functions
$u\colon \Bbb B^n \to \Bbb R^n$ that have first distributional 
derivatives $\nabla u$
such that $$\int_{\Bbb B^n} \left(|u(x)|^p+|\nabla u(x) |^p\right
)\,\,\,dx<\infty.$$ 
The $p$-capacity better suited to our problem is the 
relative first order 
variational $p$-capacity defined also in [HKM, Chapter 2]. 
We will occasionally need the Sobolev class 
$$ACL^p(\Bbb B^n)=\left\{u\in ACL(\Bbb B^n) \text{ such that } 
\int_{\Bbb B^n} |\nabla
u(x)|^p\,\,dx <\infty 
\right\}.$$ Here $ACL(\Bbb B^n)$ is the class of functions absolutely 
continuous on
almost every line. These functions are continuous and their gradients 
are Borel
functions. See for example [V\"a, \S 26]. Smooth functions 
are dense in $W^{1,p}(\Bbb B^n)$, see [K]. In particular $ACL^p(\Bbb B^n)$ 
is dense in $W^{1,p}(\Bbb B^n)$.

It was proved in [MV1] that:
\proclaim{Theorem 1.3}
 Let $u$ be  a continuous monotone function in $W^{1,p}(\Bbb B^n)$. 
 Suppose that $n-1<p\leq n$. Let $E$ be 
the set on the boundary of the unit ball 
where the non-tangential limit of $u$ does not exist, then $E$ has 
$p$-capacity zero.
\endproclaim

The proof is based on the modulus method after we obtained the following extension of Lindel\"of's theorem.
\proclaim{Theorem 1.4}
Let $u$ be  a continuous monotone function in $W^{1,p}(\Bbb B^n)$. 
Suppose that $n-1<p\leq n$. Then, for any $\epsilon>0$, there exists 
an open 
set $U$ in  $\Bbb R^n$ satisfying $cap_{p}(U)<\epsilon$ such 
that 
for any $x_0\in \partial \Bbb B^n 
\setminus U$  and $\gamma$ any curve ending at $x_0$
in $\Bbb B^n$ with $$\lim_{x\to x_0,\,\,x\in \gamma} u(x)=\alpha,$$
then $u(x)$ has  non-tangential limit $\alpha$ 
at $x_0$. 
\endproclaim  
The limitation $p>n-1$ 
appears in a module 
estimate on $(n-1)$-dimensional spheres. \par

\beginsection{ $\S$2.  Preliminaries and Oscillation Estimate.}

Let us continue with some standard notation that will be used throughout the paper. 
The open ball centered at $x_0$ with radius $r$ is denoted by 
$B^n(x_0,r)$. By $c(\alpha,\beta,\ldots)$ we denote a constant that depends only on 
the parameters $\alpha,\beta,\ldots $ and that 
may change value from line to line.
\par
Let $\Gamma$ be a family of curves in $\Bbb R^n$. Denote by $\Cal 
F(\Gamma)$ the collection
of admissible metrics for $\Gamma$. These are nonnegative Borel measurable 
functions $\rho\colon\Bbb R^n\to\Bbb R\cup\
{\lbrace \infty\rbrace}$ such that
$$\int_{\gamma}\rho\,ds\ge 1$$
for each locally rectifiable curve $\gamma\in\Gamma$. For $p\ge1$ the 
weighted $p$-module of $\Gamma$ is defined by
$$M_p(\Gamma)=\inf_{\rho\in\Cal F(\Gamma)}\int_{\Bbb 
R^n}\rho^p\,\, dx.$$
If $\Cal F(\Gamma)=\emptyset$, we set $M_p(\Gamma)=\infty$. 
Upper bounds for moduli are obtained by testing with a particular 
admissible metric.\par
Before we state our main result in this paper, let us recall some definitions.

Let $\Omega\subset\Bbb R^n$, $n\geq 2$, be a domain and 
$F\colon\Omega\to\Bbb
 R^n$ be a mapping in the Sobolev space $W_{\text{loc}}^{1,n}(\Omega; 
R^n)$ of functions in $L^n_{\text{loc}}(\Omega; R^n)$ whose distributional 
derivatives belong to $L^n_{\text{loc}}(\Omega; R^n)$. We can think of $F$ as a 
deformation of some material
whose initial configuration is $\Omega$, and we seek some functional
$I(F)$ representing the (nonlinear) elastic energy whose minimum is
attained at $F$, see [B1,2,3] and [S]. The differential of $F$ at a 
point $x$ is denoted by $DF(x)$, its norm is 
$$|DF(x)|=\sup\{|DF(x)\,h|\colon h\in \Bbb R^n,\,\,|h|=1\}$$
and its Jacobian determinant is $J_F(x)=\text{det}\,DF(x)$. We assume 
that $F$ is orientation preserving, meaning that $J_F(x)\geq 0$ for a.~e.\ 
$x \in \Omega$. The {\it dilatation} of $F$ at the point $x$ is 
defined by the ratio
$$K(x)={{|DF(x)|^n}\over{J_F(x)}}\cdot$$
If $K(x)\in L^{\infty}(\Omega; R^n),$  then $F$ is said to be a 
quasiregular mapping. \par
We will say that $F$ is a mapping of {\it finite dilatation} if
$$1\leq K(x)<\infty \text{\,\,for a.~e.~} x\in\Omega;$$
that is, except for a set of measure zero in $\Omega$, if $J_F(x)=0$ 
then $DF(x)=0$.\par
From now on we will assume that $K\in L^{n-1}$ and our domain $\Omega$ will be the 
unit ball $\Bbb B^n$ of $\Bbb R^n$.
\proclaim{Definition 1.5}
	Let $F\colon\Bbb B^n\to\Bbb
 R^n$ be a mapping. We define the multiplicity function of $F$ at some point $y\in \Bbb R^n$ 
 with respect to some domain $D\subset \Bbb B^n$ as 
 $$N(F,D,y)=\#\lbrace x\in D\colon F(x)=y\rbrace.$$
\endproclaim

Let us state our main result.
\proclaim{Theorem 1.6}
 Let $F$ be a continuous mapping in the Sobolev space $W_{\text{loc}}^{1,n}(\Bbb B^n; 
R^n)$, of finite dilatation in $W^{1,p}(\Bbb B^n)$ for 
 $n-1<p$. Let us assume also that the dilatation function $K(x) \in L^{p}$ for some $p>n-1$. 
 Let $E$ the set on the boundary of $\Bbb B^n$ where the radial limit exist and are equal to 
 $a$. Let $\Bbb B_\epsilon=\lbrace y\colon \|y\|<\epsilon\rbrace$, and $h(r)$ be a real function, a 
 constant $c$ independent of $\epsilon$ and an $\epsilon_0>0$ such that
 $$N(F,\Bbb B^n,y)\leq c h(\|y\|)$$
 for any $0<\|y\|\leq \epsilon_0.$ then if either
 $$\lim_{m\to\infty}{1\over{m^p}}\,\,\,\Biggl( \int_{e^{-(m+1)}}^{e^{-m}} {{h(r)}\over {r}}\,\,\,dr\Biggr)
 ^{{p\over n}}=0,$$
 or
 $$\sup_{\Bbb B_\epsilon} \,\,\biggl(\intav_{\Bbb B_\epsilon} N(y,\Bbb B^n,F)\,dy\biggr)\,\,\,
\biggl(\intav_{\Bbb B_\epsilon} N(y,\Bbb B^n,F)^{{1\over{1-n}}}\, dy\biggr)^{n-1}<C,$$
for any $0<\epsilon<\epsilon_0$.
 Then if $E$ has 
positive variational $p$-capacity then the mapping $F$ is identically equal to $a$.
\endproclaim
Observe that by [VG] mappings in the Sobolev space $W_{\text{loc}}^{1,n}(\Bbb B^n; 
R^n)$ of finite dilatation their component functions are monotone functions.\par
Before we pass to the proof of our result, let us examine how it is related to the results of Koskela, Mizuta and 
Miklyukov and Vuorinen. Let us observe, first, that their covering condition on $\Bbb B_\epsilon$ is an integral 
condition involving the gradient of the real function $u$ to some power $p$. For example, Koskela's condition requires 
that
$$\int_{\Bbb B_\epsilon} \|\nabla u\|^p\,\,dm\leq C\,\,\epsilon^p\,\,\biggl(\log{1\over\epsilon}\biggr)^{p-1},$$
for $1<p\leq n$. Since the function $u$ is a real function, this covering condition is for an interval $(-\epsilon, 
\epsilon)$ about zero, therefore this condition should be replaced in the case of mappings from $\Bbb B^n$ into $\Bbb 
R^n$ by an integral condition on some power of the norm of the differential matrix $\|DF\|$ or the Jacobian $J_F$.\par
As we will see later, the two extra conditions we impose on the mapping $F$ come naturally, 
first the integrability (in $L^{n-1}$) of the dilatation function $K(x)$ and second the monotonicity of the components 
of the mapping $F$. Yet, the second condition can be removed if we consider that the limits at the set 
$E\subset \partial \Bbb B^n$ are fine boundary limits. Loosely speaking, the norm of the gradient of the functions $u$ 
in the work of Koskela, Mizuta and Miklyukov and Vuorinen should be replaced by $J_F^{{1\over n}}$ in the case of a 
mapping $F$ and thus, Koskela condition on $\Bbb B_\epsilon$ translates in this case to show that (we will assume 
from now on that the mapping $F$ is sense preserving, that is $J_F\geq 0$ a.e. in $\Bbb B^n$)
$$\int_{\Bbb B_\epsilon} J_F(x)^{{p\over n}}\,\,\,dm(x)\leq C\,\,\epsilon^p\,\,\biggl(\log{1\over\epsilon}
\biggr)^{p-1}.$$
Let us quickly show here that this is the case under the hypothesis of our theorem on the multiplicity function 
$N(y,\Bbb B^n,F)$. We are trying to find an upper bound of the integral

$$\int_{\Bbb B_\epsilon} J_F(x)^{{p\over n}}\,\,\,dm(x)\leq C\,\,\,\Biggl(\int_{\Bbb B_\epsilon} J_F(x)\,\,\,dm(x)
\Biggr)^{{p\over n}},$$
where we have applied Holder's inequality. By the properties of the mapping $F$, we have that the following change of 
variable formula holds,
$$\int_{\Bbb B_\epsilon} J_F(x)\,\,\,dm(x)=\int_{y\colon \|y\|<\epsilon} N(y,\Bbb B^n,F) \,\,dm(y).$$
Therefore, by the theorem's condition on $N(y,\Bbb B^n,F)$ we have that
$$\int_{\Bbb B_\epsilon} J_F(x)\,\,\,dm(x)\leq \int_{y\colon \|y\|<\epsilon} h(r=\|y\|) \,\,dm(y),$$
and by choosing our function conveniently  i.e. $h(r)=\biggl( \log{1\over r}\biggr)^{\delta}$ for some positive $\delta$
conveniently chosen, we obtain Koskela's condition. In a similar way we can obtain Mizuta's and 
Miklyukov and Vuorinen's.\par
Mizuta's theorem in [M], states that if we have a function $\phi$ positive and nonincreasing on the interval $
(0,\infty)$ with the properties that;
$$(1)\,\,\,\,\,A^{-1}\,\,\phi(r)\leq \phi(r^2)\leq A\,\,\phi(r)$$
for every $r>0$ and $A$ a positive constant and
$$(2)\,\,\,\,\,\int_0^1 [\phi(r)]^{{{-1}\over{(p-1)}}}\,\,r^{-1}\,\,dr=\infty,$$ 
for some $1<p<\infty$. Then, the uniqueness result follows if we have a Koskela's type condition
$$(3)\,\,\,\,\,\int_{\Bbb B_\epsilon} \|\nabla u(x)\|^p\,\,dm(x)\leq C\,\,\epsilon^p\,\,\phi(\epsilon)$$
for any positive $\epsilon$. To show that Mizuta's result follows from ours we'll need to show that if we take his 
function $\phi$ to be our function $h$ and we impose his two conditions on $h$ then the above inequality follows 
from ours. Mizuta's condition (2) for $\phi$ is equivalent by Ohtsuka [O] to condition (2) in our theorem 
if we take $N(\|y\|,\Bbb B^n,F)=h(\|y\|)$, that is, Mizuta's condition (2) follows from the multiplicity function 
being an $A_n$ weight.\par
Let us pass to show that if $h$ satisfies Mizuta's conditions (1) and (2) then inequality (3) follows with 
$\|\nabla u\|$ replaced by $J_F$ when we go from real functions $u$ to mappings $F$. Thus we need to show 
that
$$ \int_{\Bbb B_\epsilon} J_F(x)^{{p\over n}}\,\,\,dm(x)\leq C\,\,\,\Biggl(\int_{\Bbb B_\epsilon} J_F(x)\,\,\,dm(x)
\Biggr)^{{p\over n}}$$
$$=\Biggl(\int_{y\colon \|y\|<\epsilon} N(y,\Bbb B^n,F) \,\,dm(y)\Biggr)^{{p\over n}}.$$
Using the fact that $N(\|y\|,\Bbb B^n,F)=h(\|y\|)$, we need to show that
$$\Biggl(\int_{y\colon \|y\|<\epsilon} h(\|y\|) \,\,dm(y)\Biggr)^{{p\over n}}\leq C\,\,\epsilon^p\,\,h(\epsilon),$$
and we are done. For this, we write the above integral as
$$\int_{y\colon \|y\|<\epsilon} h(\|y\|) \,\,dm(y)=\int_0^\epsilon  h(r) \,\,r^{n-1}dr$$
$$=\sum^\infty_{j=0} \int_{\epsilon^{2^{j+1}}}^{\epsilon^{2^{j}}}   h(r) \,\,r^{n-1}dr.$$
Now, since the function $h$ is nonincreasing and satisfies condition (1) we have that
$$\leq \sum^\infty_{j=0}  h(\epsilon^{2^{j}})\,\,{{(\epsilon^{2^{j}})^n - (\epsilon^{2^{j+1}})^n}\over n}
\leq\sum^\infty_{j=0}  h(\epsilon)\,\,A^j \,\,{{(\epsilon^{2^{j}})^n - (\epsilon^{2^{j+1}})^n}\over n}$$
$$=h(\epsilon)\,\,  \Biggl[ \epsilon^n+ (A-1)\,\,(\epsilon^n)^2\,\,\sum_{k=0}^\infty 
(A\,\,\epsilon^{2n})^k\Biggr].$$
The last series in the above equality is a geometric series whose ratio $A\,\,\epsilon^{2n}$ is less than one, 
since we have freedom in our choice of $\epsilon$, thus we have that
$$\int_{\Bbb B_\epsilon} J_F(x)^{{p\over n}}\,\,\,dm(x)\leq C\,\,\, h(\epsilon)\,\,\Biggl( \epsilon^n+
(A-1)\,\,(\epsilon^n)^2\,\,C\,\Biggr)^{{p\over n}} $$
$$\leq   C\,\,\,h(\epsilon)\,\, \Biggl(\epsilon^n\Biggr)^{{p\over n}}=C\,\,h(\epsilon)\,\,\epsilon^p,$$
as we wanted to show.\par
As for Miklyukov and Vuorinen's result, it is not difficult to show that condition (4) (for p=n) in their theorem 
is the same as condition (2) in our result (the $A_n$ weight condition for the multiplicity function of $F$ 
$N(y,\Bbb B^n,F)$, this follows from the Corollary in page 197 of Ohtsuka's paper [Oh].\par
Also, with a condition on $J_F$ as in our theorem conditions (5) and (6) in [MV] follow. We showed above how condition 
(5) follows, we will show now how condition (6) follows, that is 
$$\sum_{k=0}^\infty {1\over{h(2^{-k})}}=\infty.$$
For that, let us start as in Koskela, we assume that $\phi(r)=h(r)^{p-1}$, thus by condition (2) in [K] we have that
$$\int_0^1 {1\over {h(r)\,\,r}}\,\,dr=\infty$$
since the function $h(r)$ is nonincreasing we can rewrite the above integral as follows,
$$\infty=\int_0^1 {1\over {h(r)\,\,r}}\,\,dr=\infty=\sum_{k=0}^\infty \int_{2^{-(k+1)}}^{2^{-k}}
 {1\over {h(r)\,\,r}}\,\,dr$$
$$\leq \sum_{k=0}^\infty {1\over{h(2^{-k})}}\,\,\int_{2^{-(k+1)}}^{2^{-k}}
 {1\over {r}}\,\,dr=C\,\,\sum_{k=0}^\infty {1\over{h(2^{-k})}},$$
 and condition (6) in [MV] follows.\par
Observe also, that condition (2) in Mizuta's paper is equivalent to the necessary and sufficient condition 
in the Corollary in page 197 of Ohtsuka [Oh] once we make our choice of $\phi(r)=h(r)=N(r,\Bbb B^n,F)$.\par

After all of this, we have informally shown that our result Theorem 1.6 generalizes all the previous results 
related to the M. and F. Riesz's uniqueness theorem in two directions, namely:\par
1) Our results is for mappings.\par
\vskip 0.2in
2) When restricted to functions, i.e. $|\nabla u|$ replaced by $J_F^{{1\over n}}$ we recuperate all the known 
results by Koskela, Mizuta and Miklyukov and Vuorinen.\par

\remark{Remark 1} If we are given a function $u$ as in Koskela [K], the question will be how to extend it to a mapping 
$F$ by finding its $n$ component functions $u_i,\,\,i=1,\ldots,n$ in such a way that the conditions on the component 
functions of $F$ given by Koskela will guarantee the condition in our theorem for the mapping $F$. Obviously 
the choice of the component functios of the mapping $F$ has to be related to the function $u$. The question is, 
in which way?.
\endremark

Let us examine now the monotonicity condition on our result. In the theorem we state that 
our function ``approaches'' a fixed point $a$ on a set $E$ on the boundary $\partial \Bbb B^n$. 
We require then monotonicity, to be able to say that for any rectifiable curve $\gamma$ in $\Bbb B^n$ 
ending at a point in $E$: either the limit of the mapping $F$ along this curve does not exists 
or else, the limit exists and is equal to the same fixed point $a$. The reason we are able to conclude that 
is because as we mentioned before, for the class of mappings we are considering a Lindel\"of theorem holds, 
as can be seen in [MV1] module sets of $p$-capacity 0 on $\partial \Bbb B^n$.\par
If now, we define the term ``approaching'' as in [Z], that is, fine boundary limit exist and is equal to the fixed 
point $a$ in $\Bbb R^n$, then  we will show (it will take some work in terms of delicate estimates, that the 
same result holds without the assumption that the components of the mapping $F$ have to be monotone. We will show 
this in section 4.\par

\beginsection{ $\S$3.  Point-wise behavior of weighted Sobolev functions }

For the proof of Theorem 1.7 we are going to need the following lemma,\par
\proclaim{Lemma 2.1}
Under the hypotheses of Theorem 1.7 if we consider the weight $w_1(x)=\bigl(\ln\ln {1\over{\|F(x)\|}}\bigr)^{n-1}$, then 
there exist a positive constant $C$ independent of $p$, the point $x_0$, $w_1$ and $r$, such that
$${1\over{r^p}}\,\,\int_{\Bbb B^n(x_0,r)} w_1(x)\,\,dx\leq (p,w_1)-\text{cap} (\Bbb B^n(x_0,r);\Bbb B^n(x_0,2\,r)).$$
\endproclaim

\demo{Proof}
By a result in [HKM] all we need to show is that the measure $\mu$ defined as $d\mu(x)=w_1(x)\,dx$ satisfies 
a Poincare-type inequality. Namely, for each $\eta(x)\in C_0^\infty (\Bbb B^n(x_0,2\,r))$ we need to show that
$$\int_{\Bbb B^n(x_0,2\,r)} \eta^{n-1}(x)\,\,d\,\mu(x)\leq C\,\,r^{n-1}\,\,\int_{\Bbb B^n(x_0,2\,r)} \|\nabla 
\eta(x)\|^{n-1}\,\,d\,\mu(x),$$
where $C$ ia a constant as in the estatement of the lemma.\par
In [MV1] it was proved that the function $\eta(x)\,\,\bigl(\ln\ln {1\over{\|F(x)\|}}\bigr)$ is in the class 
$W_0^{1,n-1}\bigl(\Bbb B^n(x_0,2\,r)\bigr)$. This implies that this function satisfies a Poincare type 
inequality with respect to the euclidean metric. Thus, we have that

$$\int_{\Bbb B^n(x_0,2\,r)} \eta^{n-1}(x)\,\,\bigl(\ln\ln {1\over{\|F(x)\|}}\bigr)^{n-1}\,\,dx$$

$$\leq C\,\,r^{n-1}\,\,\int_{\Bbb B^n(x_0,2\,r)} \|\nabla 
\bigl(\eta(x)\,\,\bigl(\ln\ln {1\over{\|F(x)\|}}\bigr)\bigr)\|^{n-1}\,\,dx).$$

Applying the product rule to the right hand side of the above inequality we obtain two terms, namely

$$\int_{\Bbb B^n(x_0,2\,r)} \eta^{n-1}(x)\,\,\bigl(\ln\ln {1\over{\|F(x)\|}}\bigr)^{n-1}\,\,dx$$

$$\leq C\,\,r^{n-1}\,\,\int_{\Bbb B^n(x_0,2\,r)} \|\nabla \eta(x)\|^{n-1}
\bigl(\ln\ln {1\over{\|F(x)\|}}\bigr)\bigr)\|^{n-1}\,\,dx$$

$$+\int_{\Bbb B^n(x_0,2\,r)}  \eta(x)^{n-1}
|\nabla \bigl(\ln\ln {1\over{\|F(x)\|}}\bigr)|^{n-1}\bigr)\,\,dx.$$

By a result in [MV1] we can bound the second term inside the braces above as

$$\int_{\Bbb B^n(x_0,2\,r)}  \eta(x)^{n-1}
|\nabla \bigl(\ln\ln {1\over{\|F(x)\|}}\bigr)|^{n-1}\bigr)\,\,dx$$

$$\leq C\,\,\Biggl(
 \int_{\Bbb B^n(x_0,2\,r)} \|\nabla \eta(x)\|^{n}\,\,K(x)^{n-1}\,\,dx\Biggr)^{{{n-1}\over n}}\,\,$$

 $$\,\,\,
 \Biggl(
 \int_{\Bbb B^n(x_0,2\,r)} K(x)^{n-1}\,\,dx\Biggr)^{1\over n},$$

 since $K(x)\in L^{n-1} $ and $\eta(x)\in C_0^\infty (\Bbb B^n(x_0,2\,r))$, it easily follows that among 
 the two terms on the right hand side of the above inequality, the significant one is the first, thus we have that
 $$\int_{\Bbb B^n(x_0,2\,r)} \eta^{n-1}(x)\,\,\bigl(\ln\ln {1\over{\|F(x)\|}}\bigr)^{n-1}\,\,dx$$

$$\leq C\,\,r^{n-1}\,\,\int_{\Bbb B^n(x_0,2\,r)} \|\nabla \eta(x)\|^{n-1}
\bigl(\ln\ln {1\over{\|F(x)\|}}\bigr)\bigr)\|^{n-1}\,\,dx$$

which constitutes the desired Poincare inequality. 

 \enddemo

\proclaim{Theorem 1.7}
 Let $F$ be a continuous mapping in the Sobolev space $W_{\text{loc}}^{1,n}(\Bbb B^n; 
R^n)$, of finite dilatation in $W^{1,p}(\Bbb B^n)$. Let $\Bbb B_\epsilon=\lbrace y\colon \|y\|<\epsilon\rbrace$, 
and $h(r)$ be a real function and an $\epsilon_0>0$ such that
 $$N(F,\Bbb B^n,|y|)= h(\|y\|)$$
 for any $0<\|y\|\leq \epsilon_0.$ Then if the dilatation function $K(x) \in L^{n-1}$, 
 $F$ is discrete and open if and only if
 $$\int_0 {1\over {r\,\,h(r)^{{1\over{n-1}}}\,\,\bigl(\ln\ln {1\over{\|y\|}}\bigr)^{n}}}\,\,dr=\infty,$$
 provided that
 $$\sup_{\Bbb B_\epsilon} \,\,\biggl(\intav_{\Bbb B_\epsilon} N(|y|,\Bbb B^n,F)\,
 \bigl(\ln\ln {1\over{\|y\|}}\bigr)^{n}\,\,dy\biggr)\,$$
 $$\,\,
\biggl(\intav_{\Bbb B_\epsilon} N(|y|,\Bbb B^n,F)^{{1\over{1-n}}}\,\bigl(\ln\ln {1\over{\|y\|}}\bigr)^{{{n}\over{1-n}}} dy\biggr)^{n-1}<C,$$
for any $0<\epsilon<\epsilon_0$.
\endproclaim

\demo{Proof}
One of the implications follows from The corollary in page 197 of [Mi] and the equivalence between the weighted $p$-
modulus and the corresponding $(p,w_1)$-variational capacity where $w_1(y)=\bigl(\ln\ln {1\over{\|F(x)\|}}\bigr)^{n-1}$.\par
The other implication is more complicated to proof. Without loss of generality we can assume that $a=0$, thus we 
need to show that under our hypotheses $F^{-1}\{0\}=E$ is a discrete set, and this will follow immediately 
if we could show that $M_{n-1}^{w_1}(\Lambda(E))=0$ , since then this will imply that the 
1-Hausdorff measure of $E$ is equal to $0$ and thus the set $E$ is discrete and by a 
standard result of Titus and Young continuity of $F$ will follow. Here $\Delta=\Lambda(E)$ consists of the family 
of rectifiable curves in $\Bbb R^n$ ending at a point in $E$.\par
Let us pass to show that $M_{n-1}^{w_1}(\Lambda(E))=0$. It follows from a result in [VG] that if the mapping $F$ is 
of finite dilatation then its components are monotone and thus using a result in [MV] the class of mapping under 
consideration satisfy a Lindel\"of type theorem, which allows to say that $F(E)=\{0\}$. Consider now $\Delta_1=\Lambda\{0\}$ 
to be the set of all the rectifiable curves ending at 0. Let $\rho$ be an admissible metric for the family $\Delta_1$, 
then the metric $C\,\,(\rho\circ F)\,(x)\,\,\|DF(x)|$ is admissible for the family $\Delta$ for a suitable constant 
$C$ depending only on the dimension $n$. It is also immediate by Lindel\"of's theorem that if we have a rectifiable 
curve $\gamma\in \Delta$ then $F\circ\gamma$ belongs to the family $\Delta_1$. Let us see this with some detail.\par
Let $\gamma\in \Delta$ such that $\gamma=\lbrace (x_1(t),\ldots,x_n(t))\colon t\in [a,b]\rbrace$ then $F\circ\gamma$ 
is given in parametric equations by $\lbrace (y_1(t),\ldots,y_n(t)\colon y_i(t)=F_i(x_1(t),\ldots,x_n(t)); 
i=1,\ldots,n\rbrace$. If we denote by $s$ the arclength parameter for $\gamma$ and by $\bar s$ the arclength 
parameter for the curve $F\circ \gamma$, a simple exercise on the chain rule in several variables gives us 
that they are related by the following formula
$$d\bar s(y)\leq ds(x)\,\,\sqrt{\bigl[\sum_{i=1}^n \bigl({{\partial F_1}\over{\partial x_i}}\bigr)^2+\ldots+ 
\sum_{i=1}^n \bigl({{\partial F_n}\over{\partial x_i}}\bigr)^2\bigr]},$$
and since $\sum_{i=1}^n \bigl({{\partial F_j}\over{\partial x_i}}\bigr)^2\leq \|DF(x)\|^2$ for each $j=1,\ldots,n$ we 
obtain that
$$d\bar s(y)\leq \sqrt{n}\,\,ds(x)\,\,\|DF(x)\|.$$
Thus our constant $C=\sqrt{n}$. By the definition of the modulus of order $p$ we have that
$$M_p^{w_1}(\Delta)\leq \int C^p\,\,(\rho\circ F)^p(x)\,\,\bigl(\ln\ln {1\over{\|F(x)\|}}\bigr)^{n-1}\,\,
\|DF(x)\|^p\,\,dm(x),$$
multiplying and dividing by $K(x)^{p\over n}$ we have that
$$M_p^{w_1}(\Delta)\leq \int C^p\,\,(\rho\circ F)^p(x)\,\,\bigl(\ln\ln {1\over{\|F(x)\|}}\bigr)^{n-1}\,\,
\|DF(x)\|^p\,\,K(x)^{p\over n}\,\,
{1\over{K(x)^{p\over n}}}dm(x),$$
applying H\"older's inequality and using the fact that $J_F(x)={{\|DF(x)\|^n}\over {K(x)}}$ we obtain
$$M_p^{w_1}(\Delta)\leq \Bigg(\int C\,\,(\rho\circ F)^n(x)\,\,\bigl(\ln\ln {1\over{\|F(x)\|}}\bigr)^{{{(n-1)n}\over p}}
\,\,J_F(x)\,\,dm(x)\Bigg)^{p\over n}\,\,$$
$$\,\,\,
\Bigg(\int K(x)^{p\over {n-p}}\,\,dm(x)\Bigg)^{{{n-p}\over n}}.$$
Since $J_F(x)>0$ a.e. by assumption, we can use the formula for the change of variable to obtain that
$$M_p^{w_1}(\Delta)\leq \Bigg(\int C\,\,\rho^n(y)\,\,N(|y|,\Bbb B^n,F)\,\,
\bigl(\ln\ln {1\over{\|y\|}}\bigr)^{{{(n-1)n}\over p}}
\,\,dm(x)\Bigg)^{p\over n}\,\,$$
$$\,\,\,
\Bigg(\int K(x)^{p\over {n-p}}\,\,dm(x)\Bigg)^{{{n-p}\over n}}.$$
We now take $p=n-1$, thus
$$M_{n-1}^{w_1}(\Delta)\leq \Bigg(\int C\,\,(\rho\circ F)^n(x)\,\,N(|y|,\Bbb B^n,F)\,\,
\bigl(\ln\ln {1\over{\|F(x)\|}}\bigr)^{n}
\,\,dm(x)\Bigg)^{{{n-1}\over n}}\,\,\,\,\,
\Bigg(\int K(x)^{n-1}\,\,dm(x)\Bigg)^{{{1}\over n}},$$
taking the infimum over all the admissible metrics $rho$ for $\Delta_1$ we obtain that
$$M_p^{w_1}(\Delta)\leq  C\,\,\bigl(M_n^w(\Delta_1)\bigr)^{{{n-1}\over n}}\,\,\,\,\,
\Bigg(\int K(x)^{n-1}\,\,dm(x)\Bigg)^{{{1}\over n}},$$
where $w(y)=N(|y|,\Bbb B^n,F)\,\,\bigl(\ln\ln {1\over{\|y\|}}\bigr)^{n}$. The hypotheses of our theorem 
guaranty by the corollary in page 197 
of [Mi] that $M_n^w(\Delta_1)=0$ and since we have also assumed that $K(x)\in L^{n-1}$ this implies that 
$M_{n-1}^{w_1}(\Delta)=0$ which automatically implies that the $(n-1,w_1)$-variational capacity of $E$ is equal to 0.\par
Our aim next, will be to show that this implies that the linear measure of $E$ is equal to 0 from which the discreteness 
of $E$ will follow.\par
Let us denote by $H_1(E)$ the 1-Hausdorff measure of the set $E$. We want to show that $H_1(E)=0$. By definition,
$$H_1(E)=\lim_{\delta\to 0} \bigl[ \inf \lbrace \sum r_i\colon E\subset \bigcup \Bbb B^n(x_i,r_i),\,\,
0<r_i<\delta\rbrace\bigr].$$
Thus, applying Lemma 2.1 for $p=n-1$ and $w_1(x)=\bigl(\ln\ln {1\over{\|F(x)\|}}\bigr)^{n-1}$ we obtain that

$${1\over{r^{n-1}}}\,\,\int_{\Bbb B^n(x_0,r)} \bigl(\ln\ln {1\over{\|F(x)\|}}\bigr)^{n-1}
\,\,dx\leq (n-1,w_1)-\text{cap} (\Bbb B^n(x_0,r);\Bbb B^n(x_0,2\,r)).$$

Let $\Bbb B^n(0,\delta)$ be defined as the ball centered at 0 and radius $\delta$, and let $\Omega_\delta=
F^{-1}(\Bbb B^n(0,\delta))$. Let us take a ring completely contained in $\Omega_\delta$ centered at $x_0\in E$ and 
defined by the concentric balls $ \Bbb B^n(x_0,r)$ and $\Bbb B^n(x_0,2\,r)$ as shown in the diagram below
\vskip 2in
Observe that $|F(x)|\leq \delta$ for any $x\in \Omega_\delta$. Therefore, we have the following inequality

$${1\over{r^{n-1}}}\,\,\int_{\Bbb B^n(x_0,r)} \bigl(\ln\ln {1\over{\delta}}\bigr)^{n-1}
\,\,dx\leq {1\over{r^{n-1}}}\,\,\int_{\Bbb B^n(x_0,r)} \bigl(\ln\ln {1\over{\|F(x)\|}}\bigr)^{n-1}\,\,dx.$$

Computing the term in the left hand side of the above inequality we obtain

$$C\,\,r\,\, \bigl(\ln\ln {1\over{\delta}}\bigr)^{n-1}
\,\,\leq {1\over{r^{n-1}}}\,\,\int_{\Bbb B^n(x_0,r)} \bigl(\ln\ln {1\over{\|F(x)\|}}\bigr)^{n-1}\,\,dx,$$

where $C$ is a universal constant. Let us define the function $h(x_0,r)$ as follows

$$h(x_0,r)={1\over{r^{n-1}}}\,\,\int_{\Bbb B^n(x_0,r)} \bigl(\ln\ln {1\over{\|F(x)\|}}\bigr)^{n-1}\,\,dx.$$

With this new notation, we have that

$$C\,\,r\,\, \bigl(\ln\ln {1\over{\delta}}\bigr)^{n-1}
\,\,\leq {{h(x_0,r)}\over r}.$$

Observe that when $\delta\to 0$ we have that $x_0\to E$ and $r\to 0$. So we have that

$$\lim_{\delta\to 0} {{h(x_0,r)}\over r}=\infty.$$

That is, for any positive $\epsilon$ we can find a positive $\beta$ such that $r<\epsilon\,\,h(x_0,r)$ whenever 
$0<r<\beta$ and $x_0$ is close enough to the set $E$. Let $\lbrace \Bbb B^n(x_i,r_i)\colon x_i\in E \text{ and } 
0<r_i<\beta\rbrace_i$ be a covering of the set $E$. If we define $H_1^\delta=\inf \lbrace 
\sum r_i\colon E\subset \bigcup \Bbb B^n(x_i,r_i),\,\,
0<r_i<\delta\rbrace$, where without loss of generality we can assume that all the $x_i$'s are in $E$, we have that

$$H_1^\delta\leq \sum_i r_i\leq \sum_i \epsilon\,\, h(x_i,r_i)\leq \epsilon\,\, 
\sum_i (n-1,w_1)-\text{cap} (\Bbb B^n(x_i,r_i);\Bbb B^n(x_i,2\,r_i)).$$

We know that $(n-1,w_1)-\text{cap} (E)=0$. Hence by the definition of variational capacity and using rings to cover the 
set $E$ rather than balls (observe that we can always assume that both rings and balls are centered at points in 
$E$), then for any positive $\tilde\epsilon$ we can find a covering of $E$ by rings $E\subset 
\bigcup_i \bigl(\Bbb B^n(x_i,r_i);\Bbb B^n(x_i,2\,r_i)\bigr)$ such that 

$$\sum_i (n-1,w_1)-\text{cap} (\Bbb B^n(x_i,r_i);\Bbb B^n(x_i,2\,r_i))\leq (n-1,w_1)-\text{cap} (E) +\tilde\epsilon.$$

Choosing that covering in the previous inequality we have that

$$H_1^\delta\leq \epsilon\,\,\bigl[  (n-1,w_1)-\text{cap} (E) +\tilde\epsilon  \bigr] =\epsilon \,\,\tilde\epsilon$$
and since both $\epsilon$ and $\tilde\epsilon$ are arbitrary, letting $\delta\to 0$ we obtain that $H_1(E)=0$ as 
we wanted to show.\par
In particular, $E=F^{-1}\{0\}$ can not contain any segment, and therefore it is a totally disconnected set. 
Replacing $F(x)$ by $F(x)-b$ it follows that $F^{-1}\{b\}$ is totally disconnected for any $b\in \Bbb R^n$. The mapping 
$F$ is thus an orientation preserving light mapping and it follows from a theorem of Titus and Young [TY] that $F$ 
is open and discrete and theorem 1.7 is proved.

\qed
\enddemo

\beginsection{ $\S$4.  Proof of Theorem 1.6}

We now pass to prove Theorem 1.6.

\demo{Proof of Theorem 1.6}

Let $\rho$ be an admissible metric for $\Delta_1=\Lambda(\{0\})$, then the metric $\tilde\rho(x)=C\,\,\rho(F(x))\,\,|DF(x)|$ is 
admissible for the family $\Delta=\Lambda(E)$ for a suitable constant 
$C$ depending only on the dimension $n$, this follows as in Theorem 1.7.\par
By the definition of the modulus of order $p$ we have that

$$M_p(\Delta)\leq \int C^p\,\,(\rho\circ F)^p(x)\,\,
\|DF(x)\|^p\,\,dm(x),$$

multiplying and dividing by $K(x)^{p\over n}$ we have that

$$M_p(\Delta)\leq \int C^p\,\,(\rho\circ F)^p(x)\,\,
\|DF(x)\|^p\,\,K(x)^{p\over n}\,\,
{1\over{K(x)^{p\over n}}}dm(x),$$

applying H\"older's inequality and using the fact that $J_F(x)={{\|DF(x)\|^n}\over {K(x)}}$ we obtain

$$M_p(\Delta)\leq \Bigg(\int C\,\,(\rho\circ F)^n(x)\,\,J_F(x)\,\,dm(x)\Bigg)^{p\over n}\,\,$$

$$\,\,\,
\Bigg(\int K(x)^{p\over {n-p}}\,\,dm(x)\Bigg)^{{{n-p}\over n}}.$$

Since $J_F(x)>0$ a.e. by assumption, we can use the formula for the change of variable to obtain that

$$M_p(\Delta)\leq \Bigg(\int C\,\,\rho^n(y)\,\,N(|y|,\Bbb B^n,F)\,\,dm(x)\Bigg)^{p\over n}\,\,$$

$$\,\,\,
\Bigg(\int K(x)^{p\over {n-p}}\,\,dm(x)\Bigg)^{{{n-p}\over n}}.$$

Choose now for large positive $m$ the admissible metrics $\rho(y)=\rho_m(y)$ as follows; $\rho_m(y)={1\over{m\,\,|y|}}$ 
whenever $e^{-(m+1)}<|y|\leq e^{-m}$ and 0 otherwise. Thus we have

$$M_p(\Delta)\leq \Bigg(\int_{e^{-(m+1)}<|y|\leq e^{-m}} C\,\,{1\over{m^n\,\,|y|^n}}\,
\,N(|y|,\Bbb B^n,F)\,\,dm(x)\Bigg)^{p\over n}\,\,$$

$$\,\,\,
\Bigg(\int K(x)^{p\over {n-p}}\,\,dm(x)\Bigg)^{{{n-p}\over n}},$$

using the hypothesis of the theorem $N(|y|,\Bbb B^n,F)\leq C\,\,h(|y|)$ we have that

$$M_p(\Delta)\leq {1\over{m^p}}\,\,\Bigg(\int_{e^{-(m+1)}}^{e^{-m}} C\,\,{{h(|y|)}\over{|y|}}\,
\,dm(x)\Bigg)^{p\over n}\,\,$$

$$\,\,\,
\Bigg(\int K(x)^{p\over {n-p}}\,\,dm(x)\Bigg)^{{{n-p}\over n}}.$$

It is clear that the exponent of the dilatation function $K(x)$ in the formula above $p\over {n-p}$ is some 
$\tilde p>n-1$, so we choose our $p$ such that $p\over {n-p}=\tilde p$, hence the second factor on the right hand 
side of the above inequality is finite. Now, letting $m\to\infty$ we obtain by one of the hypothesis of our 
theorem that $M_p(\Delta)=0$. It is well known that this implies that the variational 
$p$-capacity of $E$ is equal to 0 and the theorem is proved.

\qed  
\enddemo

\beginsection{ $\S$5.  Concluding remarks}

Finally in this section, we are going to explore further the results of this paper when applied to 
some known classes of mappings. We will start with the quasiregular mappings.\par
It is well known that a mapping $F$ is quasiregular if $K(x)$ is an $L^\infty$ function, thus $1\leq K(x)\leq K<\infty$ 
a.e. x, and in this case ${1\over C}\,\,J_F(x)\leq |DF(x)|\leq C\,\, J_F(x)$ a.e. x for some constant independent of $x$ 
and we don't have to 
make any further assumption on the 
dilatation function. Our theorems 1.6 and  1.7 will then read as follows;
\proclaim{Theorem 5.1}
Let $F$ be a continuous mapping in the Sobolev space $W_{\text{loc}}^{1,n}(\Bbb B^n; 
R^n)$, of bounded dilatation. Let $\Bbb B_\epsilon=\lbrace y\colon \|y\|<\epsilon\rbrace$, 
and $h(r)$ be a real function and an $\epsilon_0>0$ such that
 $$N(F,\Bbb B^n,|y|)= h(\|y\|)$$
 for any $0<\|y\|\leq \epsilon_0.$  
 $F$ is discrete and open if and only if
 $$\int_0 {1\over {r\,\,h(r)^{{1\over{n-1}}}\,\,\bigl(\ln\ln {1\over{\|y\|}}\bigr)^{n}}}\,\,dr=\infty,$$
 provided that
 $$\sup_{\Bbb B_\epsilon} \,\,\biggl(\intav_{\Bbb B_\epsilon} N(|y|,\Bbb B^n,F)\,
 \bigl(\ln\ln {1\over{\|y\|}}\bigr)^{n}\,\,dy\biggr)\,$$
 $$\,\,
\biggl(\intav_{\Bbb B_\epsilon} N(|y|,\Bbb B^n,F)^{{1\over{1-n}}}\,\bigl(\ln\ln {1\over{\|y\|}}\bigr)^{{{n}\over{1-n}}} dy\biggr)^{n-1}<C,$$
for any $0<\epsilon<\epsilon_0$.
\endproclaim

and

\proclaim{Theorem 5.2}
 Let $F$ be a continuous mapping in the Sobolev space $W_{\text{loc}}^{1,n}(\Bbb B^n; 
R^n)$, of bounded dilatation. Let $E$ the set on the boundary of $\Bbb B^n$ where the radial limit exist 
and are equal to 
 $a$. Let $\Bbb B_\epsilon=\lbrace y\colon \|y\|<\epsilon\rbrace$, and $h(r)$ be a real function, a 
 constant $c$ independent of $\epsilon$ and an $\epsilon_0>0$ such that
 $$N(F,\Bbb B^n,y)\leq c h(\|y\|)$$
 for any $0<\|y\|\leq \epsilon_0.$ then if either
 $$\lim_{m\to\infty}{1\over{m^p}}\,\,\,\Biggl( \int_{e^{-(m+1)}}^{e^{-m}} {{h(r)}\over {r}}\,\,\,dr\Biggr)
 ^{{p\over n}}=0,$$
 or
 $$\sup_{\Bbb B_\epsilon} \,\,\biggl(\intav_{\Bbb B_\epsilon} N(y,\Bbb B^n,F)\,dy\biggr)\,\,\,
\biggl(\intav_{\Bbb B_\epsilon} N(y,\Bbb B^n,F)^{{1\over{1-n}}}\, dy\biggr)^{n-1}<C,$$
for any $0<\epsilon<\epsilon_0$.
 Then if $E$ has 
positive variational $p$-capacity then the mapping $F$ is identically equal to $a$.
\endproclaim

let us know talk about different choices of the function $h(|y|)$ which will satisfy the hypotheses of our theorems, 

which will clearly validate them. As we mentioned in previous sections, by choosing

$$h(|y|)=h(r)=\Bigl( \ln {1\over{|y|}}\Bigr)^\delta$$

for small values of $r$ and suitable choices of $\delta$, we recover all the results in [K], [Mi1] and [MV]. For 

instance choosing $\delta$ such that ${{\delta\,\,p}\over n}=p-1$ we recover the results in [K] and thus in 

[Mi1] and a straighforward computation will show that the function

$$h(|y|)=\Bigl( \ln {1\over{|y|}}\Bigr)^{{{n\,(p-1)}\over p}}$$

satisfies the conditions in Theorem 1.6.

\enddocument